\documentclass{amsproc}
\newtheorem{theorem}{Theorem}[section]

\theoremstyle{definition}
\newtheorem{definition}[theorem]{Definition}

\theoremstyle{remark}

\numberwithin{equation}{section}

\title[Reversibility]{Reversibility Questions in Groups arising in Analysis}
\author[A. O'Farrell]{Anthony G. O'Farrell}
\address{%
Mathematics Department, NUI, Maynooth, Co. Kildare, Ireland}
\email{admin@maths.nuim.ie}

\subjclass[2010]{Primary 20E34, 37-02, 46H05, 46L05}

\keywords{Reversible map, Banach algebra, Group, Dynamical system}

\thanks{The author is grateful for the 
support of the CRM, and the hospitality of Paul Gauthier
and K.N. GowriSankaran.}

\dedicatory{Dedicated to Paul Gauthier}

\begin{document}

\begin{abstract}
An element $g$ of a group is called {\em reversible}
if it is conjugate in the group to its inverse.  
In this paper we review some results
about the structure of groups involving
the reversible elements
and we pose some questions about
groups associated to a Banach algebra.
\end{abstract}

\maketitle

\section{Introduction}

\begin{definition}
An element $g$ of a group is called {\em reversible}
if it is conjugate in the group to its inverse, i.e. there
exists some map $h$, belonging to the group, such that
the conjugate $g^h=h^{-1}gh$ equals $g^{-1}$.
We say that $h$ {\em reverses} $g$, in this case.
\end{definition}


This concept had its origins in the study
of dynamical systems \cite{Bi1915}.  Classical 
conservative systems
such as the 
harmonic oscillator
and a system of 
$n$ bodies moving under their 
mutual (Newtonian) gravitational
attraction, and billiards (on a
table without pockets or corners!)
have what is called a
time-reversal symmetry: a bijection of the phase space
which conjugates the dynamical system to its inverse.

I became interested in reversibility because I encountered
reversible dynamical systems when studying a number of
different problems related to approximation:
\begin{itemize}
\item
Approximation by functions of the form
$F(x,y)=f(x)+g(y)$ on compact sets in $\mathbb R^2$ \cite{Marshall-OF1, Marshall-OF}.
\item Approximation by polynomials
of the form $p(z^2, \bar z^2+\bar z^3)$ on a disk
in the complex plane \cite{OFS}.
\end{itemize}
I was not the first to find that an apparently
undynamical problem had an essential
connection  to some discrete
dynamical system. In fact, this phenomenon has been
named {\em hidden dynamics} by V.I. Arnold
\cite{V1}.
Other problems of this type are:
\begin{itemize}
\item 
Biholomorphic classification of a pair of tangent real-analytic arcs
in the plane \cite{AhernGong, Kasner1915, Nakai, V1, V2}.
\item The polynomial hull of a disk having an isolated 
complex tangent \cite{MoWe1983, sanabria:polynomial}. 
\end{itemize}

Each of these problems involves a pair of non-commuting involutions, 
and so relates
to a reversible element in some group of maps.
Moreover, it turned out that this connection to dynamics provided
the key to resolving the problem.

With this in mind, I decided to come at the phenomenon from
the other end, and to try to understand the phenomenon
of reversibility in general.  The natural context is
the theory of groups, and I have been calculating examples,
assembling what is known from the literature, and
encouraging people to examine reversibility
in their favourite groups.  In fact, it turns out that
there is a gigantic literature on reversibility
in various different groups.  Much of this falls into
islands (or continents) of work, where the workers
are unaware of the connection between this aspect
of their field and of other fields.  The work cuts 
right across
the various fields of mainstream mathematics.
See \cite{AhernRudin, Baake-Roberts, BaRo2001, BaRo2003, 
BaRo2006, Br1996, Bullett1988, Devaney, 
Dj1967, Gomez-Meiss, GoMe2004, gong:fixed,
Go1999, Go1975, HaKa1958, Kane01, La1992, 
LaRoCa1993, Lam98, QuCa1989, Ra1981, Sarnak, Sev86, We1996, 
We1998, Wo1966}, and further references therein
and below.

In this short paper, I review some results about reversibility
and particularly factorisation
into reversible factors, and pose some questions
about groups related to Banach algebras, specifically.

I have drawn upon material in the draft of a book 
on aspects of reversibility that
is in preparation with Ian Short, and wish to acknowledge
his help with this.


\section{Notation}

Let $G$ be a group. We use the following notation:

\def\id{{\hbox{id}}}

$I=I(G):= \{f\in G: f^2=\id\}$ ---
the set of {\emph{involutions}} of $G$ (including 
the identity, $\id$).

${R_f}=R_f(G):=\{h\in G: f^h = f^{-1}\}$ (where $f^h=h^{-1}fh$) 
---the set of {\emph{reversers
of $f$}}.

$R={R(G)} := \{f\in G: R_f\not=\emptyset\}$ --- the set of {\emph{reversible
elements}}.

For $A\subset G$, ${A^n}:=\{f_1\cdots f_n:f_j\in A\}$, and  
${A^\infty} := \bigcup_{n=1}^\infty A_n$.

\medskip

Elements of $I^2$ are called {\emph{strongly-reversible.
}}They are reversed by an involution.

Membership in $I^n$ or $R^n$ is a conjugacy invariant, and $I^2\subset R$.

$I^\infty$ and $R^\infty$ are normal subgroups of $G$.

\section{Example: ${\sf GL}(n,\mathbb C)$}
Classification of linear reversible maps 
on $\mathbb C^n$ is simple.
Suppose $F\in{\sf GL}(n,\mathbb C)$ (the general linear group
over the field $\mathbb C$ of complex numbers)
 is reversible.
Since the Jordan normal form of $F^{-1}$
consists of blocks of the same size as $F$
with reciprocal eigenvalues, 
the eigenvalues of $F$ that are not $\pm1$
must split into groups of pairs $\lambda,1/\lambda$.
Furthermore, we must have the same number of Jordan blocks
of each size for $\lambda$ as for $1/\lambda$.
Vice versa, if the eigenvalues of $F$ are either $\pm1$
or split into groups of pairs $\lambda,1/\lambda$
with the same number of Jordan blocks of each size,
then both $F$ and $F^{-1}$ have the same Jordan normal form
and are therefore conjugate to each other.

Incidentally, each matrix $A\in{\sf GL}(n,\mathbb C)$ is mapped to
a conjugate of its inverse by some outer automorphism
of the group.  In fact $A\mapsto (A^t)^{-1}$ is
an automorphism, and $A^t$ is conjugate to $A$. There
is a wider context concerning \lq\lq outer reversibility"
in a group.
 
\section{The Basic Questions}
In each group, $G$, we ask:

\begin{itemize}
\item Which $f$ are reversible in $G$?

\item
Which $h$ reverse a given $f$?

\item
Describe $I^\infty$.

\item
Describe $R^\infty$.

\item
Is $I^n=I^\infty$ for some $n$?

\item
Is $R^n=R^\infty$ for some $n$? 

\item
Does every nonempty $R_g$ have
an element of finite order?  If so, what orders occur?
Is $\min\{\textup{o}(h): h\in R_g\}$
bounded, for $g\in R$?

\end{itemize}

If $g$ is reversible by some element of finite
order, then it is the product of
two elements of that (even) 
order.  Thus results about
$R^n$, combined with results
about the order of reverses,
also give information 
about factorizing elements of $G$
as a product of elements of at most 
a given order. Some find this interesting
\cite{HlOmRa2001}.

We wish to encourage people to investigate
the questions above in their favourite groups.
Many groups have been analysed, and 
in the next section we survey some 
of these results. But many
groups remain to be investigated.
In the final section, we shall draw particular
attention to groups associated to Banach algebras.

\section{Survey of Known Results}
We now give a summary of some answers, 
in examples of various categories
of groups, with a sampling of relevant sources
(by no means complete).  
Here, $G$ always denotes the group under consideration.
We give no detail about the derivation of these results,
just references to the related literature.  Some of the proofs are quite deep,
and they draw on diverse branches of mathematics. 

\subsection{}
If $G$ is abelian, then obviously 
$R=I=I^\infty$.

If $G$ is free, then it is not hard to see that
$R=I=\{1\}$.

\subsection{Finite Groups}
\cite{Co47, KolesnikovN,SiTh2008,St1998,TiZa2005}
\begin{description}
\item[\bf Dihedral  $D_n$]
$I^2=G=R$.
\item[\bf Symmetric  $S_n$]
$I^2=G=R$.
\item
[\bf Finite Coxeter]
$I^2=G=R$.
\item[\bf Alternating  $A_n$]
$I^2=R$.
\\
\quad $R\not=G$, except when $n\in\{1,2,5,6,10,14\}$.
\item[\bf Quaternion $8$-group]
$I^2\not=R=G$.
\item[\bf Finite, simple $G$]$\ $\\
\quad $R=\{1\}$ if $|G|$ is odd.
\\
\quad $G=R^2$, if $|G|$ is even, except for PSU$(3,3^2)$.
\\
\quad In general, $G\not= I^2$; when it happens is known.
\end{description}

\subsection{Classical Groups}
\cite{Ba2002, Ba1977, BuKnNi1997, Baake-Roberts, Dj1986,
DjMa1982, El1977, ElNo1982, El1993, ElMa1990, El2004b,
Go1997, HoPa1970, Ka1988, KnNi1987a, Kn1988, Lav06,
Ni1987}
\begin{description}
\item[\bf General Linear ${\sf GL}(n,F)$ ($n>1$)]
\quad $I^2=R$.
\\\quad $I^4=I^\infty = R^2$. 
\item[\bf Special Linear ${\sf SL}(n,\mathbb C)$]
\quad $I^2=R$ unless $n=2$ (mod$4$).
\\\quad $R^2=G$.
\item[\bf Orthogonal $ {\sf O}(n,\mathbb R)$ \\
($\approx$ spherical isometries)]
\quad $I^2=G$.
\item[\bf Special Orthogonal ${\sf SO}(n,\mathbb R)$]
\quad $I^2=R$.
\\\quad $I^3=G$ if $n\ge3$.
\\\quad $I^2=G$ if $n\not=2$ $(\mod4)$.
\item[\bf Unitary ${sf U}(n,\mathbb C)$]
\quad $I^2=R$
\\\quad $I^4=I^\infty$.
\item[\bf Special Unitary ${\sf SU}(n,\mathbb C)$]
\quad $I^2\not=R$.
\\\quad $I^3\not=I^6=G=R^2$.
\item[\bf Unitary Quaternionic ${\sf Sp}(n,\mathbb C)$
$={\sf Symp}(2n,\mathbb C)\cap{sf U}(2n,\mathbb C)$]$\ $\\
\quad $I^2\not=R=G=I^6$.
\item[\bf Spinor ${\sf Spin}(n,\mathbb C)$]
\quad $I^2=G$ if $n=0,1,7 \mod8$.
\\\quad $R=G$ unless $n=2 \mod4$.
\\\quad $I^4=G$ if $n\ge5$.
\end{description}

\subsection{Discrete Matrix Groups}
\cite{Baake-Roberts, BaRo2006, Is1995, La1997}

\begin{description}
\item[${\sf GL}(n,\mathbb Z)$]
\quad $I^{3n+9}=G$. $I^{41}=G$ for $n\ge84$.
\\
\quad $I^2\not=R\subset I^4$ when $n=2$.
\item[\bf Modular ${\sf PSL}(2,\mathbb Z)$]
\quad $I^2=R$.
\end{description}

\subsection{Finite-dimensional Isometry Groups}
\cite{Sh2008}
\begin{description}
\item[\bf Euclidean ${\sf Isom}(\mathbb R^n)$]
\quad $I^2=G$.
\item[\bf Orientation-preserving ${\sf Isom}^+(\mathbb R^n)$] 
\quad $I^3=G$ if $n\ge3$.
\\
\quad $I^2=G$ if $n=0$ or $3$ (mod $4$).
\item[\bf Hyperbolic ${\sf Isom}({\sf H}^n)$]
\quad $I^3=G$ if $n\ge2$.
\end{description}

\subsection{Homeomorphism Groups}
\cite{An1962, Calica, Fine-Schweigert,
FeZu1982, GiOFSh2009, Jarzcyk-b, Jarzcyk-a, OFa04,
Young1994}
\begin{description}
\item[${\sf Homeo}(\mathbb R)$]
\quad $I^2\not=R$.
\quad $I^3\not=I^4=G=R^2$.
\item[${\sf Homeo}^+(\mathbb R)$]
\quad $R^4=G$.
\item[${\sf Homeo}(\mathbb S^1)$]
\quad $I^2\not=R$.
\quad $I^3=R^2=G$.
\item[${\sf Homeo}^+(\mathbb S^1)$]
\quad $I^2\not=R$.
\quad $I^3=R^2=G$.
\item[${\sf Homeo}(\mathbb S^n)$]
\quad $G=I^6$ when $n=2$ or $3$. (Open for $n>3$).
\item[\bf Compact surface MCG]
\quad $I^n\not=G=I^\infty$, $\forall n\in\mathbb N$,
if genus $>2$.
\end{description}

\subsection{Maps with extra Structure}
\cite{AhernOF, AR, Birkhoff-1939, 
Banakh-Yagasaki, GiSh2010, OFShort}
\begin{description}
\item[\bf Diffeomorphism ${\sf Diffeo}(\mathbb R)$]
\quad $I^2\not=R$.
\\
\quad $I^3\not=I^4=G=R^2$.
\item[${\sf Diffeo}^+(\mathbb R)$]
\quad $R^4=G$.
\item[\bf Formal germs on $(\mathbb C^n,0)$]\cite{OFa08, AhernOF}
\quad $I^4=I^\infty$ when $n=1$.
\\
\quad $R^2=R^\infty$ when $n=1$.
\\
\quad $R^k=R^\infty$ with $k= 3+2\cdot\textup{ceiling}(\log_2 n)$
\\
\quad \quad when $n\ge2$.
\\
\quad $I^{15}=R^\infty$ when $n=2$.
\item[\bf Piecewise-linear ${\sf PL}(\mathbb R)$]\cite{BrSq01, BrSq1985}
\quad $I^2\not=R$.
\\
\quad $I^3\not=I^4=G=R^2$.
\item[${\sf PL}^+(\mathbb R)$]
\quad $I=\{1\}$.
\quad $R^4=G$.
\item[\bf PL with finitely many nodes ${\sf PLF}(\mathbb R)$]
\quad $I^2=R$.
\item[${\sf PLF}^+(\mathbb R)$]
\quad $R^4=G$.
\end{description}

\section{Banach Algebras}

Let $A$ be a Banach algebra.
We may associate two collections of groups to $A$.
In each case, we ask the usual questions. 

\subsection{$A^{-1}$ and its subgroups}

Suppose $A$ has identity (or adjoin one, if not) and $\|1\|=1$.

Reversibility in $A^{-1}$ is not interesting unless
$A$ is noncommutative. Also central reversibles are
just central involutions, so the real problems
are about the quotient
$$ \frac{A^{-1}}{Z(A^{-1})} \equiv {\sf Inn}(A). $$

One interesting subgroup is
$${\sf Iso}(A) = \{x\in A: \|x\| = \|x^{-1}\|=1\}.$$
This coincides with the subgroup (often denoted ${sf U}(A)$ \cite{Allan}) 
of unitary elements, in case $A$
is a $C^*$ algebra.

One may also focus on the connected component
of $1$ in either group, $G$, and on the
intersection $G\cap E^c$ with
the commutator of any subset $E\subset A$.

One also has the 
normal subgroup \\
$\{x\in A: \|a-1\|<1\}^\infty$,
which lies in the group $(\exp A)^\infty$.

For instance, Gustafson, Halmos, and Radjavi
\cite{GuHaRa1976} showed that for finite-dimensional 
(real or complex) Hilbert spaces $H$,
the group $G={\sf GL}(H)$ has
$I^4= I^\infty$, and they noted
that for infinite-dimensional
$H$, we have $I^4\not=I^7={\sf GL}(H)$.
I don't know whether $7$ is
the best possible value in that statement.
What happens with other $C^*$ algebras? 

\medskip
It is known that for finite-dimensional 
${\sf GL}(H)$, we have $I^2=R$ in 
$G$, and also in the unitary subgroup.
What about other $C^*$ algebras?

\subsection{Aut$(A)$ and its subgroups}

The main interesting subgroup
(apart from the inner automorphism group,
already mentioned)
is the group of isometric isomorphisms.
This is often the same as Aut$(A)$.

\medskip
As an example, when $X$ is a locally-compact Hausdorff space
and $A=C_0(X,\mathbb C)$, then Aut$(A)$ is isomorphic to
${\sf Homeo}(X)$, so we have seen answers in case
$X=\mathbb R^1$ and $X=\mathbb S^1$.  
For the disk algebra, the automorphism group
is isomorphic to ${\sf PSL}(2,\mathbb R)$, so all
elements are reversible.  Similarly, for the
polydisk algebra, the automorphism group
is isomorphic to the group of conformal
automorphisms of the polydisk, which is
isomorphic to the wreath product 
$S_n\wr {\sf PSL}(2,\mathbb R)$, and again all elements
are reversible.  

\medskip
We remark that the algebra of all formal power
series in $n$ indeterminates, with complex
coefficients, has a natural Frechet algebra
structure, and embeds in some Banach algebras \cite{Allan}.
Its automorphism group is isomorphic to the 
group of formally-invertible formal
germs mentioned already.


\begin{thebibliography}{20}


\bibitem[AG05]{AhernGong}
P.~Ahern and X.~Gong, \textit{A complete classification for pairs of real
  analytic curves in the complex plane with tangential intersection}, J. Dyn.
  Control Syst. \textbf{11} (2005), no.~1, 1--71.

\bibitem[All11]{Allan}
G.R. Allan, \textit{Introduction to Banach algebras}, CUP, 2011.

\bibitem[And62]{An1962}
R.~D. Anderson, \textit{On homeomorphisms as products of conjugates of a given
  homeomorphism and its inverse}, Topology of 3-manifolds and related topics
  ({P}roc. {T}he {U}niv. of {G}eorgia {I}nstitute, 1961), Prentice-Hall,
  Englewood Cliffs, N.J., 1962, pp.~231--234.

\bibitem[AO09]{AhernOF}
P.~Ahern and A.~G. O'Farrell, \textit{Reversible biholomorphic germs}, Comput.
  Methods Funct. Theory \textbf{9} (2009), no.~2, 473--484.

\bibitem[AR95a]{AR}
P.~Ahern and J-P. Rosay, \textit{Entire functions, in the classification of
  differentiable germs tangent to the identity, in one or two variables},
  Trans. Amer. Math. Soc. \textbf{347} (1995), no.~2, 543--572.

\bibitem[AR95b]{AhernRudin}
P.~Ahern and W.~Rudin, \textit{Periodic automorphisms of {${\bf C}\sp n$}},
  Indiana Univ. Math. J. \textbf{44} (1995), no.~1, 287--303.

\bibitem[Bak02]{Ba2002}
A.~Baker, \textit{Matrix groups}, Springer Undergraduate Mathematics Series,
  Springer-Verlag London Ltd., London, 2002, An introduction to Lie group
  theory.

\bibitem[Bal78]{Ba1977}
C.~S. Ballantine, \textit{Products of involutory matrices. {I}}, Linear and
  Multilinear Algebra \textbf{5} (1977/78), no.~1, 53--62.

\bibitem[Bir15]{Bi1915}
G.~D. Birkhoff, \textit{The restricted problem of three bodies}, Rend. Circ. Mat.
  Palermo \textbf{39} (1915), 265--334.

\bibitem[Bir39]{Birkhoff-1939}
\bysame, \textit{D\'eformations analytiques et fonctions auto-\'equivalentes},
  Ann. Inst. H. Poincar\'e \textbf{9} (1939), 51--122.

\bibitem[BKN97]{BuKnNi1997}
F.~B{\"u}nger, F.~Kn{\"u}ppel, and K.~Nielsen, \textit{Products of symmetries in
  unitary groups}, Linear Algebra Appl. \textbf{260} (1997), 9--42.

\bibitem[BR97]{Baake-Roberts}
M.~Baake and J.~A.~G. Roberts, \textit{Reversing symmetry group of {${\rm
  Gl}(2,{\bf Z})$} and {${\rm PGl}(2,{\bf Z})$} matrices with connections to
  cat maps and trace maps}, J. Phys. A \textbf{30} (1997), no.~5, 1549--1573.

\bibitem[BR01]{BaRo2001}
\bysame, \textit{Symmetries and reversing symmetries of toral automorphisms},
  Nonlinearity \textbf{14} (2001), no.~4, R1--R24.

\bibitem[BR03]{BaRo2003}
\bysame, \textit{Symmetries and reversing symmetries of area-preserving
  polynomial mappings in generalised standard form}, Phys. A \textbf{317}
  (2003), no.~1-2, 95--112.

\bibitem[BR06]{BaRo2006}
\bysame, \textit{The structure of reversing symmetry groups}, Bull. Austral.
  Math. Soc. \textbf{73} (2006), no.~3, 445--459.

\bibitem[Bri96]{Br1996}
M.~G. Brin, \textit{The chameleon groups of {R}ichard {J}. {T}hompson:
  automorphisms and dynamics}, Inst. Hautes \'Etudes Sci. Publ. Math. (1996),
  no.~84, 5--33 (1997).

\bibitem[BS85]{BrSq1985}
M.~G. Brin and C.~C. Squier, \textit{Groups of piecewise linear homeomorphisms of
  the real line}, Invent. Math. \textbf{79} (1985), no.~3, 485--498.

\bibitem[BS01]{BrSq01}
\bysame, \textit{Presentations, conjugacy, roots, and centralizers in groups of
  piecewise linear homeomorphisms of the real line}, Comm. Algebra \textbf{29}
  (2001), no.~10, 4557--4596.

\bibitem[Bul88]{Bullett1988}
S.~Bullett, \textit{Dynamics of quadratic correspondences}, Nonlinearity
  \textbf{1} (1988), no.~1, 27--50.

\bibitem[BY09]{Banakh-Yagasaki}
T.~Banakh and T.~Yagasaki, \textit{The diffeomorphism groups of the real line are
  pairwise bihomeomorphic}, Topology \textbf{48} (2009), no.~2-4, 119--129.

\bibitem[Cal71]{Calica}
A.~B. Calica, \textit{Reversible homeomorphisms of the real line}, Pacific J.
  Math. \textbf{39} (1971), 79--87.

\bibitem[Cox47]{Co47}
H.~S.~M. Coxeter, \textit{The product of three reflections}, Quart. Appl. Math.
  \textbf{5} (1947), 217--222.

\bibitem[Dev76]{Devaney}
R.~L. Devaney, \textit{Reversible diffeomorphisms and flows}, Trans. Amer. Math.
  Soc. \textbf{218} (1976), 89--113.

\bibitem[Djo67]{Dj1967}
D.~{\v{Z}}. Djokovi{\'c}, \textit{Product of two involutions}, Arch. Math.
  (Basel) \textbf{18} (1967), 582--584.

\bibitem[Djo86]{Dj1986}
\bysame, \textit{Pairs of involutions in the general linear group}, J. Algebra
  \textbf{100} (1986), no.~1, 214--223.

\bibitem[DM82]{DjMa1982}
D.~{\v{Z}}. Djokovi{\'c} and J.~G. Malzan, \textit{Products of reflections in
  {${\rm U}(p,\,q)$}}, Mem. Amer. Math. Soc. \textbf{37} (1982), no.~259,
  vi+82.

\bibitem[Ell77]{El1977}
Erich~W. Ellers, \textit{Bireflectionality in classical groups}, Canad. J. Math.
  \textbf{29} (1977), no.~6, 1157--1162.

\bibitem[Ell93]{El1993}
E.~W. Ellers, \textit{The reflection length of a transformation in the unitary
  group over a finite field}, Linear and Multilinear Algebra \textbf{35}
  (1993), no.~1, 11--35.

\bibitem[EM90]{ElMa1990}
E.~W. Ellers and J.~Malzan, \textit{Products of reflections in the kernel of the
  spinorial norm}, Geom. Dedicata \textbf{36} (1990), no.~2-3, 279--285.

\bibitem[EN82]{ElNo1982}
E.~W. Ellers and W.~Nolte, \textit{Bireflectionality of orthogonal and symplectic
  groups}, Arch. Math. (Basel) \textbf{39} (1982), no.~2, 113--118.

\bibitem[EV04]{El2004b}
E.~W. Ellers and O.~Villa, \textit{The special orthogonal group is
  trireflectional}, Arch. Math. (Basel) \textbf{82} (2004), no.~2, 122--127.

\bibitem[FS55]{Fine-Schweigert}
N.~J. Fine and G.~E. Schweigert, \textit{On the group of homeomorphisms of an
  arc}, Ann. of Math. (2) \textbf{62} (1955), 237--253.

\bibitem[FZ82]{FeZu1982}
W.~Feit and G.~J. Zuckerman, \textit{Reality properties of conjugacy classes in
  spin groups and symplectic groups}, Algebraists' homage: papers in ring
  theory and related topics ({N}ew {H}aven, {C}onn., 1981), Contemp. Math.,
  vol.~13, Amer. Math. Soc., Providence, R.I., 1982, pp.~239--253.

\bibitem[GHR76]{GuHaRa1976}
W.~H. Gustafson, P.~R. Halmos, and H.~Radjavi, \textit{Products of involutions},
  Linear Algebra and Appl. \textbf{13} (1976), no.~1/2, 157--162, Collection of
  articles dedicated to Olga Taussky Todd.

\bibitem[GM03]{Gomez-Meiss}
A.~G{\'o}mez and J.~D. Meiss, \textit{Reversible polynomial automorphisms of the
  plane: the involutory case}, Phys. Lett. A \textbf{312} (2003), no.~1-2,
  49--58.

\bibitem[GM04]{GoMe2004}
\bysame, \textit{Reversors and symmetries for polynomial automorphisms of the
  complex plane}, Nonlinearity \textbf{17} (2004), no.~3, 975--1000.

\bibitem[Gon96]{gong:fixed}
X.~Gong, \textit{Fixed points of elliptic reversible transformations with
  integrals}, Ergodic Theory Dynam. Systems \textbf{16} (1996), no.~4,
  683--702.

\bibitem[Goo97]{Go1997}
G.~R. Goodson, \textit{The inverse-similarity problem for real orthogonal
  matrices}, Amer. Math. Monthly \textbf{104} (1997), no.~3, 223--230.

\bibitem[Goo99]{Go1999}
\bysame, \textit{Inverse conjugacies and reversing symmetry groups}, Amer. Math.
  Monthly \textbf{106} (1999), no.~1, 19--26.

\bibitem[GOS09]{GiOFSh2009}
N.~Gill, A.~G. O'Farrell, and I.~Short, \textit{Reversibility in the group of
  homeomorphisms of the circle}, Bull. Lond. Math. Soc. \textbf{41} (2009),
  no.~5, 885--897.

\bibitem[Gow75]{Go1975}
R.~Gow, \textit{Real-valued characters of solvable groups}, Bull. London Math.
  Soc. \textbf{7} (1975), 132.

\bibitem[GS10]{GiSh2010}
N.~Gill and I.~Short, \textit{Reversible maps and composites of involutions in
  groups of piecewise linear homeomorphisms of the real line}, Aequationes
  Math. \textbf{79} (2010), no.~1-2, 23--37.

\bibitem[HK58]{HaKa1958}
P.~R. Halmos and S.~Kakutani, \textit{Products of symmetries}, Bull. Amer. Math.
  Soc. \textbf{64} (1958), 77--78.

\bibitem[HOR01]{HlOmRa2001}
M.~Hladnik, M.~Omladi{\v{c}}, and H.~Radjavi, \textit{Products of roots of the
  identity}, Proc. Amer. Math. Soc. \textbf{129} (2001), no.~2, 459--465.

\bibitem[HP71]{HoPa1970}
F.~Hoffman and E.~C. Paige, \textit{Products of two involutions in the general
  linear group}, Indiana Univ. Math. J. \textbf{20} (1970/1971), 1017--1020.

\bibitem[Ish95]{Is1995}
H.~Ishibashi, \textit{Involutary expressions for elements in {${\rm GL}\sb n({\bf
  Z})$} and {${\rm SL}\sb n({\bf Z})$}}, Linear Algebra Appl. \textbf{219}
  (1995), 165--177.

\bibitem[Jar02a]{Jarzcyk-b}
W.~Jarczyk, \textit{Reversibility of interval homeomorphisms without fixed
  points}, Aequationes Math. \textbf{63} (2002), no.~1-2, 66--75.

\bibitem[Jar02b]{Jarzcyk-a}
\bysame, \textit{Reversible interval homeomorphisms}, J. Math. Anal. Appl.
  \textbf{272} (2002), no.~2, 473--479.

\bibitem[Kan01]{Kane01}
R.~Kane, \textit{Reflection groups and invariant theory}, CMS Books in
  Mathematics/Ouvrages de Math\'ematiques de la SMC, 5, Springer-Verlag, New
  York, 2001.

\bibitem[Kas15]{Kasner1915}
E.~Kasner, \textit{Conformal classification of analytic arcs or elements:
  {P}oincar\'e's local problem of conformal geometry}, Trans. Amer. Math. Soc.
  \textbf{16} (1915), no.~3, 333--349.

\bibitem[KN87]{KnNi1987a}
F.~Kn{\"u}ppel and K.~Nielsen, \textit{On products of two involutions in the
  orthogonal group of a vector space}, Linear Algebra Appl. \textbf{94} (1987),
  209--216.

\bibitem[KN05]{KolesnikovN}
S.~G. Kolesnikov and Ja.~N. Nuzhin, \textit{On strong reality of finite simple
  groups}, Acta Appl. Math. \textbf{85} (2005), no.~1-3, 195--203.

\bibitem[Kn{\"u}88]{Kn1988}
F.~Kn{\"u}ppel, \textit{Products of involutions in orthogonal groups},
  Combinatorics '86 ({T}rento, 1986), Ann. Discrete Math., vol.~37,
  North-Holland, Amsterdam, 1988, pp.~231--247.

\bibitem[Laf97]{La1997}
T.~J. Laffey, \textit{Lectures on integer matrices}, Unpublished lecture notes
  (1997).

\bibitem[Lam92]{La1992}
J.~S.~W. Lamb, \textit{Reversing symmetries in dynamical systems}, J. Phys. A
  \textbf{25} (1992), no.~4, 925--937.

\bibitem[Liu88]{Ka1988}
K.~M. Liu, \textit{Decomposition of matrices into three involutions}, Linear
  Algebra Appl. \textbf{111} (1988), 1--24.

\bibitem[LOS07]{Lav06}
R.~L{\'a}vi{\v{c}}ka, A.~G. O'Farrell, and I.~Short, \textit{Reversible maps in
  the group of quaternionic {M}\"obius transformations}, Math. Proc. Cambridge
  Philos. Soc. \textbf{143} (2007), no.~1, 57--69.

\bibitem[LR98]{Lam98}
J.~S.~W. Lamb and J.~A.~G. Roberts, \textit{Time-reversal symmetry in dynamical
  systems: a survey}, Phys. D \textbf{112} (1998), no.~1-2, 1--39,
  Time-reversal symmetry in dynamical systems (Coventry, 1996).

\bibitem[LRC93]{LaRoCa1993}
J.~S.~W. Lamb, J.~A.~G. Roberts, and H.~W. Capel, \textit{Conditions for local
  (reversing) symmetries in dynamical systems}, Phys. A \textbf{197} (1993),
  no.~3, 379--422.

\bibitem[MO79]{Marshall-OF1}
D.~E. Marshall and A.~G. O'Farrell, \textit{Uniform approximation by real
  functions}, Fund. Math. \textbf{54} (1979), 203--11.

\bibitem[MO83]{Marshall-OF}
\bysame, \textit{Approximation by a sum of two algebras. {T}he lightning bolt
  principle}, J. Funct. Anal. \textbf{52} (1983), no.~3, 353--368.

\bibitem[MW83]{MoWe1983}
J.~K. Moser and S.~M. Webster, \textit{Normal forms for real surfaces in {${\bf
  C}\sp{2}$} near complex tangents and hyperbolic surface transformations},
  Acta Math. \textbf{150} (1983), no.~3-4, 255--296.

\bibitem[Nak98]{Nakai}
I.~Nakai, \textit{The classification of curvilinear angles in the complex plane
  and the groups of {$\pm$} holomorphic diffeomorphisms}, Ann. Fac. Sci.
  Toulouse Math. (6) \textbf{7} (1998), no.~2, 313--334.

\bibitem[Nie87]{Ni1987}
K.~Nielsen, \textit{On bireflectionality and trireflectionality of orthogonal
  groups}, Linear Algebra Appl. \textbf{94} (1987), 197--208.

\bibitem[O'F04]{OFa04}
A.~G. O'Farrell, \textit{Conjugacy, involutions, and reversibility for real
  homeomorphisms}, Irish Math. Soc. Bull. (2004), no.~54, 41--52.

\bibitem[O'F08]{OFa08}
\bysame, \textit{Composition of involutive power series, and reversible series},
  Comput. Methods Funct. Theory \textbf{8} (2008), no.~1-2, 173--193.

\bibitem[OS09]{OFShort}
A.~G. O'Farrell and I.~Short, \textit{Reversibility in the diffeomorphism group
  of the real line}, Publ. Mat. \textbf{53} (2009), no.~2, 401--415.

\bibitem[OSG02]{OFS}
A.G O'Farrell and A.~Sanabria-Garcia, \textit{De Paepe's disc has nontrivial
  polynomial hull}, Bull. LMS \textbf{34} (2002), 490--4.

\bibitem[QC89]{QuCa1989}
G.~R.~W. Quispel and H.~W. Capel, \textit{Local reversibility in dynamical
  systems}, Phys. Lett. A \textbf{142} (1989), no.~2-3, 112--116.

\bibitem[Rad81]{Ra1981}
H.~Radjavi, \textit{The group generated by involutions}, Proc. Roy. Irish Acad.
  Sect. A \textbf{81} (1981), no.~1, 9--12.

\bibitem[Sar07]{Sarnak}
P.~Sarnak, \textit{Reciprocal geodesics}, Analytic number theory, Clay Math.
  Proc., vol.~7, Amer. Math. Soc., Providence, RI, 2007, pp.~217--237.

\bibitem[Sev86]{Sev86}
M.~B. Sevryuk, \textit{Reversible systems}, Lecture Notes in Mathematics, vol.
  1211, Springer-Verlag, Berlin, 1986.

\bibitem[SG00]{sanabria:polynomial}
A.~Sanabria-Garc{\'{\i}}a, \textit{Polynomial hulls of smooth discs: a survey},
  Irish Math. Soc. Bull. (2000), no.~45, 135--153.

\bibitem[Sho08]{Sh2008}
I.~Short, \textit{Reversible maps in isometry groups of spherical, {E}uclidean
  and hyperbolic space}, Math. Proc. R. Ir. Acad. \textbf{108} (2008), no.~1,
  33--46.

\bibitem[ST08]{SiTh2008}
A.~Singh and M.~Thakur, \textit{Reality properties of conjugacy classes in
  algebraic groups}, Israel J. Math. \textbf{165} (2008), 1--27.

\bibitem[Ste98]{St1998}
A.~Stein, \textit{{$1\frac12$}-generation of finite simple groups}, Beitr\"age
  Algebra Geom. \textbf{39} (1998), no.~2, 349--358.

\bibitem[TZ05]{TiZa2005}
P.~H. Tiep and A.~E. Zalesski, \textit{Real conjugacy classes in algebraic groups
  and finite groups of {L}ie type}, J. Group Theory \textbf{8} (2005), no.~3,
  291--315.

\bibitem[Vor81]{V1}
S.~M. Voronin, \textit{Analytic classification of germs of conformal mappings
  {$({\bf C},\,0)\rightarrow ({\bf C},\,0)$}}, Funktsional. Anal. i Prilozhen.
  \textbf{15} (1981), no.~1, 1--17, 96.

\bibitem[Vor82]{V2}
\bysame, \textit{Analytic classification of pairs of involutions and its
  applications}, Funktsional. Anal. i Prilozhen. \textbf{16} (1982), no.~2,
  21--29, 96.

\bibitem[Web96]{We1996}
S.~M. Webster, \textit{Double valued reflection in the complex plane}, Enseign.
  Math. (2) \textbf{42} (1996), no.~1-2, 25--48.

\bibitem[Web98]{We1998}
\bysame, \textit{Real ellipsoids and double valued reflection in complex space},
  Amer. J. Math. \textbf{120} (1998), no.~4, 757--809.

\bibitem[Won66]{Wo1966}
M.~J. Wonenburger, \textit{Transformations which are products of two
  involutions}, J. Math. Mech. \textbf{16} (1966), 327--338.

\bibitem[You94]{Young1994}
S.~W. Young, \textit{The representation of homeomorphisms on the interval as
  finite compositions of involutions}, Proc. Amer. Math. Soc. \textbf{121}
  (1994), no.~2, 605--610.


\end{thebibliography}
\end{document}